\documentclass[11pt,twoside]{article}

\usepackage{amsbsy,amsfonts,amsmath,amssymb,eucal,mathrsfs}
\usepackage[all]{xy}
\usepackage{pstricks,epsfig}
\usepackage{pst-plot}
\usepackage{psfrag}

\newtheorem{definition}{Definition}[section]

\newtheorem{prop}[definition]{Proposition}

\newtheorem{lemm}[definition]{Lemma}
\newtheorem{fact}[definition]{Fact}

\newtheorem{coro}[definition]{Corollary}
\newtheorem{theo}[definition]{Theorem}
\newtheorem{notation}[definition]{Notation}

\newtheorem{construction}[definition]{Construction}

\newtheorem{remark}[definition]{Remark}
\newenvironment{rema}{\begin{remark} \rm}{\end{remark}}
\newtheorem{remarks}[definition]{Remarks}

\newtheorem{example}[definition]{Example}

\newtheorem{examples}[definition]{Examples}

\newtheorem{nothing}[definition]{$\!\!$}

\newenvironment{proo}{{\flushleft \it Proof.}}{\hfill $\square$
  \vspace{2mm}}
\newenvironment{proo-prop}{{\flushleft \it Proof of Proposition
    \ref{prop-f_4}.}}{\hfill $\square$ \vspace{2mm}}
\newenvironment{conj}{\begin{conjecture} \rm}{\end{conjecture}}
\newtheorem{conjecture}[definition]{Conjecture}
\newenvironment{stat}{\begin{statement} \rm}{\end{statement}}
\newtheorem{statement}[definition]{Statement}

\newtheorem{definition*}{Definition}[section]
\newenvironment{defi*}{\begin{definition*} \rm}{\end{definition*}}
\newtheorem{definitions*}[definition*]{Definitions}
\newenvironment{defis*}{\begin{definitions*} \rm}{\end{definitions*}}
\newtheorem{prop*}[definition*]{Proposition}
\newtheorem{lemm*}[definition*]{Lemma}
\newtheorem{coro*}[definition*]{Corollary}
\newtheorem{theo*}[definition*]{Theorem}
\newtheorem{remark*}[definition*]{Remark}
\newenvironment{rema*}{\begin{remark*} \rm}{\end{remark*}}
\newtheorem{remarks*}[definition*]{Remarks}
\newenvironment{remas*}{\begin{remarks*} \rm}{\end{remarks*}}
\newtheorem{example*}[definition*]{Example}
\newenvironment{exam*}{\begin{example*} \rm}{\end{example*}}
\newtheorem{examples*}[definition*]{Examples}
\newenvironment{exams*}{\begin{examples*} \rm}{\end{examples*}}
\newtheorem{nothing*}[definition*]{$\!\!$}
\newenvironment{noth*}{\begin{nothing*} \rm}{\end{nothing*}}

\newtheorem{commentaire*}[definition*]{Commentaire}

\begin{document}

\def \JK {{\mathfrak{J}}}
\def \IK {{\mathfrak{I}}}
\def \Rh {{\widehat{R}}}
\def \Sh {{\widehat{S}}}
\def \supp {{\rm Supp}}
\def \codim {{\rm codim}}
\def \b {{\beta}}
\def \T {{\Theta}}
\def \t {{\theta}}
\def \L {{\cal L}}
\def \sca #1#2{\langle #1,#2 \rangle}
\def\pt{\{{\rm pt}\}}
\def\x {{\underline{x}}}
\def\y {{\underline{y}}}
\def\aut{{\rm Aut}}
\def\ra{\rightarrow}
\def\s{\sigma}\def\OO{\mathbb O}\def\PP{\mathbb P}\def\QQ{\mathbb Q}
 \def\CC{\mathbb C} \def\ZZ{\mathbb Z}\def\JO{{\mathcal J}_3(\OO)}
\newcommand{\G}{\mathbb{G}}
\def\proof{\noindent {\it Proof.}\;}
\def\qed{\hfill $\square$}
\def \uh {{\widehat{u}}}
\def \vh {{\widehat{v}}}
\def \fh {{\widehat{f}}}
\def \wh {{\widehat{w}}}
\def \Wh {{{W_{{\rm aff}}}}}
\def \Wt {{\widetilde{W}_{{\rm aff}}}}
\def \Qt {{\widetilde{Q}}}
\def \Waff {{W_{{\rm aff}}}}
\def \Waffm {{W_{{\rm aff}}^-}}
\def \Wpaff {{{(W^P)}_{{\rm aff}}}}
\def \Wtpaff {{{(\widetilde{W}^P)}_{{\rm aff}}}}
\def \Wtaffm {{\widetilde{W}_{{\rm aff}}^-}}
\def \lh {{\widehat{\lambda}}}
\def \pit {{\widetilde{\pi}}}
\def \lt {{{\lambda}}}
\def \xh {{\widehat{x}}}
\def \yh {{\widehat{y}}}
\def \a {\alpha}
\def \b {\beta}
\def \l {\lambda}
\def \t {\theta}
%%%%%ICI ATTENTION \def \T {\theta}

%%%%%%%%% Commandes pe %%%%%%%%%%%%%%

\newcommand{\expxy}{\exp_{x \to y}}
\newcommand{\drat}{d_{\rm rat}}
\newcommand{\dmax}{d_{\rm max}}
\newcommand{\Dmax}{D_{\rm max}}
\newcommand{\zl}{Z(x,L_x,y,L_y)}

% alg{\~A}{\AA}{!`}bres norm{\~A}{\copyright}es et anneaux usuels

\newcommand{\tiff}{if and only if }
\newcommand{\N}{\mathbb{N}}
\newcommand{\A}{{\mathbb{A}_{\rm Aff}}}
\newcommand{\Ah}{{\mathbb{A}_{\rm Aff}}}
\newcommand{\At}{{\widetilde{\mathbb{A}}_{\rm Aff}}}
\newcommand{\Ht}{{{H}^T_*(\Omega K^{\ad})}}
\renewcommand{\H}{{\rm Hi}}
\newcommand{\Ih}{{I_{\rm Aff}}}
\newcommand{\psit}{{\widetilde{\psi}}}
\newcommand{\xit}{{\widetilde{\xi}}}
\newcommand{\Jt}{{\widetilde{J}}}
\newcommand{\Zt}{{\widetilde{Z}}}
\newcommand{\Xt}{{\widetilde{X}}}
\newcommand{\at}{{\widetilde{A}}}
\newcommand{\Z}{\mathbb Z}
\newcommand{\R}{\mathbb{R}}
\newcommand{\Q}{\mathbb{Q}}
\newcommand{\C}{\mathbb{C}}
\renewcommand{\O}{\mathbb{O}}
\newcommand{\F}{\mathbb{F}}
\newcommand{\p}{\mathbb{P}}
\newcommand{\co}{{\cal O}}
\newcommand{\pos}{{\bf P}}

\renewcommand{\a}{{\alpha}}
\newcommand{\az}{\a_\Z}
\newcommand{\ak}{\a_k}

\newcommand{\rc}{\R_\C}
\newcommand{\cc}{\C_\C}
\newcommand{\hc}{\H_\C}
\newcommand{\oc}{\O_\C}

\newcommand{\rk}{\R_k}
\newcommand{\ck}{\C_k}
\newcommand{\hk}{\H_k}
\newcommand{\ok}{\O_k}

\newcommand{\rz}{\R_Z}
\newcommand{\cz}{\C_Z}
\newcommand{\hz}{\H_Z}
\newcommand{\oz}{\O_Z}

\newcommand{\RR}{\R_R}
\newcommand{\CR}{\C_R}
\newcommand{\HR}{\H_R}
\newcommand{\OR}{\O_R}

\newcommand{\re}{\mathtt{Re}}

\newcommand{\matttr}[9]{
\left (
\begin{array}{ccc}
{} \hspace{-.2cm} #1 & {} \hspace{-.2cm} #2 & {} \hspace{-.2cm} #3 \\
{} \hspace{-.2cm} #4 & {} \hspace{-.2cm} #5 & {} \hspace{-.2cm} #6 \\
{} \hspace{-.2cm} #7 & {} \hspace{-.2cm} #8 & {} \hspace{-.2cm} #9
\end{array}
\hspace{-.15cm}
\right )   }

% alg{\~A}{\AA}{!`}bre

\newcommand{\dual}{{\bf v}}
\newcommand{\com}{\mathtt{Com}}
\newcommand{\rg}{\mathtt{rg}}
\newcommand{\pu}{{\mathbb{P}^1}}
\newcommand{\scal}[1]{\langle #1 \rangle}
\newcommand{\MK}[2]{{\overline{{\rm M}}_{#1}(#2)}}
\newcommand{\mor}[2]{{{\rm Mor}_{#1}(\pu,#2)}}

\newcommand{\fg}{\mathfrak g}
\newcommand{\fgad}{{\mathfrak g}^{\rm ad}}
\renewcommand{\fh}{\mathfrak h}
\newcommand{\fu}{\mathfrak u}
\newcommand{\fz}{\mathfrak z}
\newcommand{\fn}{\mathfrak n}
\newcommand{\fe}{\mathfrak e}
\newcommand{\fp}{\mathfrak p}
\newcommand{\ft}{\mathfrak t}
\newcommand{\fl}{\mathfrak l}
\newcommand{\fq}{\mathfrak q}
\newcommand{\fsl}{\mathfrak {sl}}
\newcommand{\fgl}{\mathfrak {gl}}
\newcommand{\fso}{\mathfrak {so}}
\newcommand{\fsp}{\mathfrak {sp}}
\newcommand{\ff}{\mathfrak {f}}

\newcommand{\ad}{{\rm ad}}
\newcommand{\jad}{{j^\ad}}
\newcommand{\id}{{\rm id}}

%%%% Poids et racines   %%%%

\newcommand{\dynkinadeux}[2]
{
$#1$
\setlength{\unitlength}{1.2pt}
\hspace{-3mm}
\begin{picture}(12,3)
\put(0,3){\line(1,0){10}}
\end{picture}
\hspace{-2.4mm}
$#2$
}

\newcommand{\mdynkinadeux}[2]
{
\mbox{\dynkinadeux{#1}{#2}}
}

\newcommand{\dynkingdeux}[2]
{
$#1$
\setlength{\unitlength}{1.2pt}
\hspace{-3mm}
\begin{picture}(12,3)
\put(1,.8){$<$}
\multiput(0,1.5)(0,1.5){3}{\line(1,0){10}}
\end{picture}
\hspace{-2.4mm}
$#2$
}

\newcommand{\poidsesix}[6]
{
\hspace{-.12cm}
\left (
\begin{array}{ccccc}
{} \hspace{-.2cm} #1 & {} \hspace{-.3cm} #2 & {} \hspace{-.3cm} #3 &
{} \hspace{-.3cm} #4 & {} \hspace{-.3cm} #5 \vspace{-.13cm}\\
\hspace{-.2cm} & \hspace{-.3cm} & {} \hspace{-.3cm} #6 &
{} \hspace{-.3cm} & {} \hspace{-.3cm}
\end{array}
\hspace{-.2cm}
\right )      }

\newcommand{\copoidsesix}[6]{
\hspace{-.12cm}
\left |
\begin{array}{ccccc}
{} \hspace{-.2cm} #1 & {} \hspace{-.3cm} #2 & {} \hspace{-.3cm} #3 &
{} \hspace{-.3cm} #4 & {} \hspace{-.3cm} #5 \vspace{-.13cm}\\
\hspace{-.2cm} & \hspace{-.3cm} & {} \hspace{-.3cm} #6 &
{} \hspace{-.3cm} & {} \hspace{-.3cm}
\end{array}
\hspace{-.2cm}
\right |      }

\newcommand{\poidsesept}[7]{
\hspace{-.12cm}
\left (
\begin{array}{cccccc}
{} \hspace{-.2cm} #1 & {} \hspace{-.3cm} #2 & {} \hspace{-.3cm} #3 &
{} \hspace{-.3cm} #4 & {} \hspace{-.3cm} #5 & {} \hspace{-.3cm} #6
\vspace{-.13cm}\\
\hspace{-.2cm} & \hspace{-.3cm} & {} \hspace{-.3cm} #7 &
{} \hspace{-.3cm} & {} \hspace{-.3cm}
\end{array}
\hspace{-.2cm}
\right )      }

\newcommand{\copoidsesept}[7]{
\hspace{-.12cm}
\left |
\begin{array}{cccccc}
{} \hspace{-.2cm} #1 & {} \hspace{-.3cm} #2 & {} \hspace{-.3cm} #3 &
{} \hspace{-.3cm} #4 & {} \hspace{-.3cm} #5 & {} \hspace{-.3cm} #6
\vspace{-.13cm}\\
\hspace{-.2cm} & \hspace{-.3cm} & {} \hspace{-.3cm} #7 &
{} \hspace{-.3cm} & {} \hspace{-.3cm}
\end{array}
\hspace{-.2cm}
\right |      }

\newcommand{\poidsehuit}[8]{
\hspace{-.12cm}
\left (
\begin{array}{cccccc}
{} \hspace{-.2cm} #1 & {} \hspace{-.3cm} #2 & {} \hspace{-.3cm} #3 &
{} \hspace{-.3cm} #4 & {} \hspace{-.3cm} #5 & {} \hspace{-.3cm} #6 &
{} \hspace{-.3cm} #7   \vspace{-.13cm}\\
\hspace{-.2cm} & \hspace{-.3cm} & {} \hspace{-.3cm} #8 &
{} \hspace{-.3cm} & {} \hspace{-.3cm}
\end{array}
\hspace{-.2cm}
\right )      }

\newcommand{\copoidsehuit}[8]{
\hspace{-.12cm}
\left |
\begin{array}{cccccc}
{} \hspace{-.2cm} #1 & {} \hspace{-.3cm} #2 & {} \hspace{-.3cm} #3 &
{} \hspace{-.3cm} #4 & {} \hspace{-.3cm} #5 & {} \hspace{-.3cm} #6 &
{} \hspace{-.3cm} #7  \vspace{-.13cm}\\
\hspace{-.2cm} & \hspace{-.3cm} & {} \hspace{-.3cm} #8 &
{} \hspace{-.3cm} & {} \hspace{-.3cm}
\end{array}
\hspace{-.2cm}
\right |      }

\newcommand{\im}{\mathtt{Im}}

%%%% Caligraphic letters %%%%%%%%%%%%%%

\def\cA{{\cal A}} \def\cC{{\cal C}} \def\cD{{\cal D}} \def\cE{{\cal E}}
\def\cF{{\cal F}} \def\cG{{\cal G}} \def\cH{{\cal H}} \def\cI{{\cal I}}
\def\cK{{\cal K}} \def\cL{{\cal L}} \def\cM{{\cal M}} \def\cN{{\cal N}}
\def\cO{{\cal O}}
\def\cP{{\cal P}} \def\cQ{{\cal Q}} \def\cT{{\cal T}} \def\cU{{\cal U}}
\def\cV{{\cal V}} \def\cX{{\cal X}} \def\cY{{\cal Y}} \def\cZ{{\cal Z}}

\def \g {{\gamma}}

\def \pp {\odot}
\def \tr {{}^t}
\def \ct {{}_Tc}
\def \lcom {$\Lambda$-(co)minuscule }
\def \lmin {$\Lambda$-minuscule }
\renewcommand{\ok}[2]{{#1} \cdot {#2} = {#1} \pp {#2}}
\newcommand{\notok}[2]{{#1} \cdot {#2} \not = {#1} \pp {#2}}
\newcommand{\oknu}[3]{c_{{#1},{#2}}^{{#3}} = t_{{#1},{#2}}^{{#3}} \cdot
m_{{#1},{#2}}^{{#3}}}
\newcommand{\notoknu}[3]{c_{{#1},{#2}}^{{#3}} \not = t_{{#1},{#2}}^{{#3}} \cdot
m_{{#1},{#2}}^{{#3}}}

 \title{Rationality of some
Gromov-Witten varieties\\
and application to quantum $K$-theory}
 \author{P.-E. Chaput, N. Perrin}

\maketitle

\begin{abstract}

We show that for any minuscule or cominuscule homogeneous space
$X$, the Gromov-Witten
varieties of degree $d$ curves passing through three general points of $X$
are rational or empty for any $d$.
Applying techniques of A. Buch and L. Mihalcea
\cite{BM} to constructions of the authors together with L. Manivel \cite{CMP},
we deduce that the equivariant $K$-theoretic three points Gromov-Witten
invariants are equal to classical equivariant $K$-theoretic invariants on
auxilliary spaces again homogeneous under ${\rm Aut}(X)$.

\end{abstract}

 {\def\thefootnote{\relax}
 \footnote{ \hspace{-6.8mm}
 }
 }

\begin{center}{\bf Introduction}\end{center}

To compute Gromov-Witten invariants in 
cominuscule Grassmannians of classical type, %s $A,B,C,D$, ICI
A. Buch, A. Kresch, and
H. Tamvakis
proved a comparison formula which they called ``quantum to classical
principle'' \cite{Buch,BKT}. This formula compares Gromov-Witten invariants
on a homogeneous space $G/P$ with classical intersection numbers on
an auxilliary homogeneous space under the same group $G$.
In \cite{CMP}, L. Manivel and us gave a uniform treatment of their
results
in order to include the two exceptional homogeneous spaces
$E_6/P_1$ and $E_7/P_7$. This ``quantum to classical principle'' was further
generalised by A. Buch and L. Mihalcea
to the $K$-theoretic invariants in \cite{BM}, in the case of
Grassmannians. In this paper we show such a principle for all cominuscule
homogeneous spaces: see Corollary \ref{coro-principe}.

As it is explained in \cite{BM}, such a result is a consequence of the fact
that the locus $GW_d(x_1,x_2,x_3)$ of curves of degree $d$ passing through
three generic points $x_1,x_2,x_3$ is rational.
Denoting $\dmax$ the minimal integer $d$ such that
there passes a curve of degree $d$ through any two points in $X$,
they %therefore ICI 
made the following conjecture:

\begin{conj}({\bf Buch-Mihalcea})
  Let $X$ be a cominuscule variety and let $d > \dmax$.
  If $x_1,x_2,x_3$ are general points in
  $X$, then the Gromov-Witten variety $GW_d(x_1,x_2,x_3)$ is
  rational.
\end{conj}

This conjecture is true in type $A$ (already if $d \geq \dmax$) by
\cite[Corollary 2.2]{BM}, and we show that it is true for any cominuscule
space
except exactly in the case where $X=E_6/P_1$ and
$d=3$. To cover this case we introduce the following notation: let
$X^{(3)}_d \subset X^3$
denote the subvariety of triples of points $(x_1,x_2,x_3)$
such that there exists a stable curve of degree $d$ through
$x_1,x_2,x_3$. Proposition \ref{prop-stat-true} contains in particular:

\begin{theo}
\label{theo-rationnel}
  Let $X$ be a cominuscule variety, and $d$ be any integer.
  Then for $(x_1,x_2,x_3)$ generic in $X^{(3)}_d$, the Gromov-Witten variety
  $GW_d(x_1,x_2,x_3)$ is rational.
\end{theo}

As a consequence, we obtain a ``quantum to classical'' principle for 
the equivariant quantum $K$-theory as follows. We define (see
Section \ref{setion-k-quant-class}) for each degree $d$ a space $Y_d$
homogeneous under the group ${\rm Aut}(X)$ and consider the
incidence diagram
$$\xymatrix{Z_d\ar[r]^p\ar[d]_q&X\\
Y_d.&\\}$$
Denoting by $I_d^T(\a,\beta,\gamma)$ the equivariant $K$-theoretic
Gromov-Witten invariant and using the general framework developped
in \cite{BM}, we obtain the following:
\begin{theo}
For $X$ a cominuscule homogeneous space and $d$ any non negative
integer, we have the following formula:
$$I_d^T(\a,\beta,\gamma) 
= \chi_{Y_d}(q_*p^*\a\cdot q_*p^*\beta\cdot q_*p^*\gamma) $$
where $\chi_{M_d}:K_T(M_d)\to K_T({\rm pt})$ is the $T$-equivariant Euler 
characteristic and where $\a$, $\beta$ and $\gamma$ are three classes in 
$K_T(X)$.
\end{theo}

\vskip .3cm

We now explain our techniques to get rationality.
By \cite[Proposition 3.17, Fact 3.18]{CMP}, for $d \leq \dmax$,
this Gromov-Witten variety is one point. On the other hand, the cases of type
$A$ have already been proved in \cite{BM}. Thus we stick to the other cases.

We have two different arguments, each of them being not enough alone to cover
all the
cases. On the one hand, we explain that the quadrics, the symplectic
Grassmannians of type $C$, the spinor varieties of even type $D_{2n}$, and
the exceptional case $E_7/P_7$, share some common geometric properties
which imply that the restriction to $X$ of the projection with center an
osculating space is a birational morphism to a
projective space, and that we can describe the set of curves in this projective
space which are the projection of some curve in $X$ (Proposition
\ref{prop-courbe}). From this follows an
explicit birational
description of the Gromov-Witten varieties, and in particular the
fact that they are rational.

On the other hand, if the Dynkin diagram of $X$ has a cominuscule node
which is not the marked node (which occurs in all cases but the odd
quadrics, the
symplectic Grassmannians and $E_7/P_7$), then we can apply results of
\cite{courbes-rat} to show that a big open subset of $X$ is the total space
of a globally generated vector bundle $E$ on a smaller homogeneous space $Y$.
It follows that a dense open subset of the space of rational curves
of degree $d$ in $X$ is
a vector bundle over the space of rational curves of degree $d$ in $Y$.
We then show that on a
generic rational curve of degree $d$ in $Y$ the restriction of $E$
is the vector
bundle $\oplus \cO_{\p^1}(a_i)$ with $a_i \geq 2$ for all $i$,
which implies the rationality
of the Gromov-Witten locus by induction.

At the end of Section \ref{subsection-passing}, we also explain how the
above vector bundle technique allows us to prove some rationality
results for Gromov-Witten varieties for non cominuscule homogeneous
spaces and large degrees (this question was also raised in
\cite{BM}). We deal with orthogonal Grassmannians and
adjoint varieties (see Proposition \ref{prop-ortho} and Proposition
\ref{prop-adj}).

\tableofcontents

\def \S {{\mathbb S}}
\def \ifftext {if and only if }

\section{The case of generalised Veronese curves}
\label{section-veronese}

In this section we prove that Gromov-Witten varieties are rational for a class
of homogeneous spaces that share common geometric properties and that we
call, after Mukai \cite{mukai}, generalised Veronese curves. These homogeneous
spaces are: the quadrics, the Lagrangian Grassmannians
$\G_\omega(n,2n)$, the Grassmannians
$\G(n,2n)$, the ``even'' spinor varieties $\G_Q(2n,4n)$, and $E_7/P_7$.

More precisely, let $a \geq 1$ and $n \geq 2$ be
integers. Assume that $a \in \{1,2,4,8\}$ if $n=3$, and that
$a \in \{1,2,4\}$ if $n \geq 4$. We define a variety $X(a,n)$ by the following:

$$
\begin{array}{r|c|c}
n=2 & a \mbox{ arbitrary} & \Q^a \\
\hline 
n \geq 3 & a=1 & \G_\omega(n,2n) \\
n \geq 3 & a=2 & \G(n,2n) \\
n \geq 3 & a=4 & \G_Q(2n,4n) \\
\hline
n=3 & a=8 & E_7/P_7
\end{array}
$$

\begin{rema}
Let $X$ be a cominuscule homogeneous space (see \cite{CMP} for
more about cominuscule homogeneous spaces). Denote $\dmax$
the smallest integer such that through any two points of $X$ there
passes a rational curve of degree $\dmax$.
The varieties $X(a,n)$ are all the cominuscule homogeneous spaces $X$ such that
through any {\em three} points passes a rational curve of degree $\dmax$.
\end{rema}

\subsection{Generalised Veronese curves}

To prove that $X(n,a)$ has rational Gromov-Witten loci, we give a
parametrisation of the set of curves in $X(n,a)$ of a given degree passing
through a point $x \in X(n,a)$. To this end our first task is to show that the
projection from an osculating space to $x$ at order $n-2$ is a birational
map to
a projective space, by interpreting $X(n,a)$ as a 
generalised rational normal curve of
degree $n$.

\vskip .5cm

More precisely, let $G(n,a)$ be the fundamental covering
of the automorphism group of $X(n,a)$, $V(n,a)$ the representation of
$G(n,a)$ where $X(n,a)$ is minimally embedded,
$L(n,a)$ the derived group of
a Levi factor of the parabolic subgroup of $G(n,a)$
stabilising a point in $X(n,a)$, and $T(n,a)$ a tangent space of $X(n,a)$,
which is a representation of $L(n,a)$. We display the above defined groups 
and representations in the following array.
$$
\begin{array}{c|c|c|c|c|c}
n & a & G(n,a) & V(n,a) & L(n,a) & T(n,a) \\
\hline
2 & \mbox{arbitrary} & Spin_{a+2} & \C^{a+2} & Spin_a & \C^a \\
\hline
n \geq 3 & 1 & Sp_{2n} & (\wedge^n \C^{2n})_\omega & SL_n & S^2\C^n \\
n \geq 3 & 2 & SL_{2n} & \wedge^n \C^{2n} & SL_n \times SL_n & \C^n \otimes 
\C^n \\
n \geq 3 & 4 & Spin_{4n} & \S_+ & SL_{2n} & \wedge^2 \C^{2n} \\
\hline
3 & 8 & E_7 & V^{56} & E_6 & V^{27}
\end{array}
$$
($\S^+$ denotes a spinor representation of $Spin_{4n}$)

\vskip .2cm
\noindent
Moreover let $Z(n,a)$ be the center of a Levi factor of a parabolic subgroup
of $G(n,a)$ stabilising a point in $X(n,a)$: we have $Z(n,a) \simeq \C^*$. 
This group induces a decomposition
$$V(n,a) = V_{-n} \oplus V_{-n+2} \oplus \cdots \oplus V_{n-2}\oplus V_n$$ 
as a $L(n,a)$-representation, where $V_i$ is the $i$-th eigenspace of $Z(n,a)$.
For example, for $a=8$ we have
$V^{56} = \C \oplus V^{27} \oplus V^{27} \oplus \C$, which is known as the
Zorn-matrix expression for elements in $V^{56}$ (see for example 
\cite[Example 2.3 and Theorem 4.15]{garibaldi}). 
The other cases can be checked directly. The following fact is an easy 
exercice on plethysm:

\begin{fact}
\label{fait-puissance}
For $2 \leq i \leq n$
the space of $L(n,a)$-equivariant maps of
degree $i$ from $T(n,a)$ to $V_{-n+2i}$ is isomorphic to $\C$.
\end{fact}

The choice of a specific such map does not really matter, but to fix the
ideas we denote $\nu_i$ one degree $i$ map from $T(n,a)$ to $V_{-n+2i}$.
In the case $n \geq 3$ and $a=2$, we have $T(n,a) = \C^n \otimes \C^n$ and
$V_{-n+2i} = \wedge^i \C^n \otimes \wedge^i \C^n$, thus we denote $\nu_i(x)$
the element given by the $i\times i$ minors of the matrix
$x \in \C^n \otimes \C^n$. By restriction to symmetric matrices we also
define $\nu_i$ in case $a=1$. For skew-symmetric matrices (case $a=4$),
we define $\nu_i(x)$ for $x \in \wedge^2 \C^{2n}$ as the collection of
the Pfaffians of the symmetric $2i\times 2i$ minors of $x$.

The interpretation of $X(n,a)$ as a generalised Veronese curve of degree $n$
comes from the following proposition which is proved in 
\cite[Section 3.1]{manivel}:

\begin{prop}
\label{prop-veronese}
Let $v \in V(n,a)$ and denote $v=(v_k)$ according to the above decomposition
$V(n,a) = \oplus_k V_k$. Assume that $v_{-n} = 1$, then the class of $v$ in
$\p V(n,a)$ belongs to $X(n,a)$ \ifftext $v_{-n+2i} = \nu_i(v_{-n+2})$ for
$2 \leq i \leq n$.
\end{prop}
\begin{proo}
For self-completeness we include a quick elementary proof of this result which
is part of the general construction performed in \cite{manivel}.
When $n=2$ this result is obvious, and when $n=3$ and $a=8$ it is proved in
\cite{mukai}. Thus we assume $n \geq 3$ and $a \in \{1,2,4\}$.
Let us first consider the case when $a=2$.
Consider a base $(e_1,\ldots,e_n,f_1,\ldots,f_n)$ of $\C^{2n}$ and a
completely decomposable form $v \in \wedge^n \C^{2n}$, \emph{i.e.} 
the class of $v$ in $\p \wedge^n \C^{2n}$
is an element of the Grassmannian $\G(n,2n)$, and write $v=(v_k)$ in the 
decomposition $V=\oplus_kV_k$. Denoting $E$ resp. $F$ the subspace of 
$\C^{2n}$ generated by the $e_i$'s resp. the
$f_i$'s, we have $\C^{2n} = E \oplus F$ and 
$V_{-n+2i} = \wedge^{n-i}E \otimes \wedge^i F$. Since by assumption
$v_{-n} \not = 0$, the coordinate of $v \in \wedge^n \C^{2n}$ on
$e_1 \wedge \cdots \wedge e_n$ does not vanish, and thus $v$
can be written as
$ ( e_1 + \varphi(e_1) ) \wedge \cdots \wedge ( e_n + \varphi(e_n) ) $,
where $\varphi : E \to F$ is a linear map. The element $v_{-n+2}$ is then
$$\sum_i e_1 \wedge \cdots \wedge \varphi(e_i) \wedge \cdots \wedge e_n;$$
this element in $\wedge^{n-1} E \otimes F$
identifies with $\varphi \in {\rm Hom}(E,F)$. The element $v_{-n+4}$ is given
by 
$$
\sum_{i < j} e_1 \wedge \cdots \wedge \varphi(e_i) \wedge \cdots
\wedge \varphi(e_j) \wedge \cdots \wedge e_n:
$$
it corresponds to the element in ${\rm Hom}(\wedge^2E,\wedge^2F)$ given by the
2 by 2 minors of $\varphi$, namely $\nu_2(v_{-n+2})$.
Similarly we have $v_{-n+2i} = \nu_i(v_{-n+2})$
for any $i$. Thus the proposition is proved in the case of Grassmannians.

The case $a=1$ follows using the injections $S^2 \C^n \subset \C^n
 \otimes \C^n$
and $X(n,1) \subset X(n,2)$. The case $a=4$ also follows from the injections
$\wedge^2 \C^{2n} \subset \C^{2n} \otimes \C^{2n}$ and
$\G_Q(2n,4n) \subset \G(2n,4n)$, recalling that the Pl{\"u}cker embedding
$\G_Q(2n,4n) \subset \p \wedge^{2n} \C^{4n}$ is the second Veronese
embedding of $\G_Q(2n,4n) \subset \p V(n,a)$, and that the squares of the
Pfaffians of the symmetric submatrices are the minors.
\end{proo}

\subsection{Projection from osculating spaces}

Let $x \in \p V(n,a)$ resp. $y \in \p V(n,a)$
be the projectivisation of the highest resp. lowest weight
line, these are points in $X(n,a)$. Moreover we can assume that $x$
resp. $y$ is the class in $\p V(n,a)$ of
$(1,0,\ldots,0)$ resp. $(0,\ldots,0,1)$ in the decomposition 
$V(n,a) = \oplus_k V_k$.
From Proposition \ref{prop-veronese}
we get immediately the following corollary:

\begin{coro}
\label{coro-projection}
We have $T^{n-2}_xX = \bigoplus_{i\leq n-2} V_{-n+2i}$ and
$T^1_yX = V_{n-2} \oplus V_n$. Moreover, the restriction to $X$
of the projection to $\p T^1_yX$ with
center $\p T^{n-2}_xX$ is birational.
\end{coro} 

We denote $Y \subset \p V_{n-2} \subset \p T^1_yX$ the closed $L(n,a)$-orbit
in $\p V_{n-2}$. It has dimension $(n-1)a$.
This corollary enables one to caracterise the curves in $\p T^1_yX$
which are projections of curves in $X$:

\begin{prop}
\label{prop-courbe}
A rational curve of degree $d-n+1$ in $\p T^1_yX$
is the projection of a rational
curve of degree $d$ in $X$ passing through $x$ \ifftext it passes through
$d-n$ points in $Y \subset \p T^1_yX$.
\end{prop}
\begin{proo}
We denote $p : \p V(n,a) \dasharrow \p T^1_yX$ the projection of Corollary
\ref{coro-projection}.
Let $C \subset X$ be a rational curve of degree $d$ passing through $x$;
we show that $p(C)$ is a rational curve of degree $d-n+1$ which
passes through $d-n$ points in $Y$. In fact, since $x \in C \subset X$ and
$T^{n-2}_xX$ is the $(n-2)$-th osculating space, the intersection
$C \cap T^{n-2}_xX$ has multiplicity at least $n-1$, and generically exactly
$n-1$, hence $p(C)$ has degree $d-n+1$.

On the other hand, the intersection of $C$ with the hyperplane $T^{n-1}_xX$
has length $d$, and the multiplicity of $x$ is $n$; thus there are
$d-n$ other intersection points of $C$ with $T^{n-1}_xX$. Let $z$ be any of
these points and write $z=(z_i)$ with $z_k \in V_k$. If $C$ is generic
then $z_{-n} \not = 0$, so that we can assume that $z_{-n} = 1$, and
Proposition \ref{prop-veronese} implies that $\nu_n(z_{-n+2}) = z_n = 0$
and $z_{n-2} = \nu_{n-1}(z_{-n+2})$. Thus it follows that
$[z_{n-2}] \in Y$. Thus for $C$ a generic curve $p(C)$ meets $Y$ along
at least $d-n$ points, and therefore this holds for any $C$.

Observe now that the space of curves of degree $d-n+1$ passing
through $d-n$ points
in $Y$ is irreducible. In fact, we can specify such a curve by choosing any
$d-n$ points on $Y$ and then a curve of degree $d-n+1$ in $\p T^1_yX$
passing through these points. Since the space of curves of degree $d-n+1$
through $d-n$ fixed points in a projective space is non empty and irreducible,
the claim is proved.

To finish the proof of the
proposition, it is enough to show that a generic curve of degree
$d-n+1$ passing through $d-n$ points in $Y$ is the projection of a curve
of degree $d$ in $X$. To this end a dimension count argument is enough.
Note that we have $c_1(X)=an+2-a$
and $\dim(X)=\frac{1}{2}(an^2+2n-an)$. Thus the space of curves of degree $d$
in $X$ and passing through $x$ has dimension
$d(an+2-a)-2$. On the other hand the space of
curves of degree $d-n+1$ in $\p(T_y^1)$ passing through $d-n$ points in 
$Y$ has dimension
$(d-n+1)(\dim T^1_yX)+\dim \p T^1_yX-3-(d-n)\codim_{\p T^1_yX}(Y)$. This 
is equal
to $(d-n+2)(\dim X+1)-4-(d-n)\left ( n-1 + \frac{1}{2}(n-1)(n-2)a \right )$.
This is also
$d(an+2-a)-2$. Thus the proposition is proved.
\end{proo}

\begin{theo}
\label{theo-veronese}
Let $X=X(n,a)$ be a generalised Veronese curve and $d$ an integer. The space 
of rational curves of degree $d$ passing through 3 fixed points in $X$ is 
empty for $d<n$ and rational for $d\geq n$.
\end{theo}
\begin{proo}
For $d<n$ an easy dimension count proves the result. 

For $d\geq n$, by Proposition \ref{prop-courbe}, a generic curve of degree 
of $d$ passing through 3 fixed points in $X$ corresponds to a generic curve 
of degree $d-n+1$ passing through 2 fixed points in $\p T^1_yX$ and passing 
moreover through $d-n$ points in $Y$. Thus it is enough to show that this 
latter space is rational.

This in turn follows from the fact $Y$ is rational and
that in any projective space the space
of rational curves of degree $\delta$ ($\delta=d-n+1$) passing through
$\delta+1$ fixed points is non empty and rational.
\end{proo}

\def \GW {{\rm GW}}
\newcommand{\mattt}[9]{
\left (
\begin{array}{rcl}
{} #1 & {} #2 & {} #3 \\
{} #4 & {} #5 & {} #6 \\
{} #7 & {} #8 & {} #9 \\
\end{array}
\right )
}

\section{Rationality by the vector bundle technique}
\label{subsection-passing}

\subsection{Inductive birational realisation of homogeneous varieties}
\label{section-fibre}

In this subsection, we recall some results of \cite{courbes-rat} that we shall
use in the sequel. Let $Q$ be a parabolic subgroup of $G$ containing $B$. Let
$w_0$ be the longest element of the Weyl group. Consider the subvariety
$Z=Q^{w_0}/Q^{w_0}\cap P$ of $G/P$ where $Q^{w_0}$ denotes the conjugate of
$Q$ under $w_0$. Because $Q^{w_0}$ contains the opposite Borel subgroup $B^-$,
the subvariety $Z$ is open in $X$. Furthermore, denoting by $i$ the Weyl
involution, one easily checks (for more details see \cite[Propposition
6]{courbes-rat}) that the complement of $Z$ has codimension at least 2 as soon
as the sets of simple roots $\Sigma(P)$ and $i(\Sigma(Q))$ are disjoint.

Let us consider the Levi decomposition of $Q$ given by $Q=U_Q\rtimes L_Q$ where
$U_Q$ is the unipotent radical of $Q$ and $L_Q$ is a section of the quotient
$Q/U_Q$. This decomposition induces a decomposition of $Q^{w_0}$ as follows
$Q^{w_0}=U_Q^{w_0}\rtimes L_Q^{w_0}$. Let us denote by $Y$ the quotient
$L_Q^{w_0}/L_Q^{w_0}\cap P$, we have a natural morphism $p:Z\to Y$ induced by
the projection of $Q^{w_0}$ to $L^{w_0}_Q$. As a consequence of
\cite[Proposition 5]{courbes-rat} we have the

\begin{prop}
\label{prop-fibre}
  Let $Q$ be a maximal parabolic subgroup such that $\Sigma(Q)$ is a
  cominuscule root, then the map $p$ realises $Z$ as the vector bundle over
  $L_Q^{w_0}/L_Q^{w_0}\cap P$ defined by the representation of $L_Q^{w_0}\cap
  P$ in $U_Q^{w_0}/U_Q^{w_0}\cap P$.
\end{prop}

\begin{proo}
  This statement is contained in \cite[Proposition 5]{courbes-rat} with the
  exception that $p$ is in general an affine bundle. However, as explained at
  the end of the proof of \cite[Proposition 5]{courbes-rat}, this affine
  bundle is an actual vector bundle as soon as the unipotent subgroup
  $U_Q^{w_0}$ is abelian. This is the case as soon as $\Sigma(Q)$ is a
  cominuscule root.
\end{proo}

\begin{rema}
\label{rema-global-gen}
  As explained in \cite[Proposition 5]{courbes-rat}, the above vector bundle
  $p:Z\to Y$ is globally generated.
\end{rema}

\subsection{Passing rationality of Gromov-Witten varieties through vector
  bundles}

In this subsection, we let $p:X\to Y$ be the vector bundle defined in the
previous subsection and we explain how to deduce the rationality of
Gromov-Witten varieties in $Z$ from the rationality of Gromov-Witten varieties
in $Y$ under simple hypothesis on the vector bundle. Let us denote by $E$ the
locally free sheaf on $Y$ defined by taking local sections of the vector
bundle $p:Z\to Y$. Let us assume that the Picard group of $Y$ and hence of $Z$
is $\Z$.

\begin{prop}
\label{prop-rat-fibre}
  Let $x$, $y$ and $z$ be three points in general position in $Z$ and assume
  that $d$ is a non negative integer such that there exists a degree $d$
  morphism $f:\pu\to Y$ with $H^1(\pu,f^*E\otimes\cO_\pu(-3))=0$. 

Assume that the Gromov-Witten variety of degree $d$ rational curves in $Y$
passing through $p(x)$, $p(y)$ and $p(z)$ is rational, then the Gromov-Witten
variety of degree $d$ rational curves in $Z$ passing through $x$, $y$ and $z$
is also rational.
\end{prop}

\begin{proo}
  The ideas for the proof are partially already contained in the proof of
  \cite[Proposition 3]{courbes-rat}. Indeed, the map $p$ induces by
    composition a morphism $q:\mor{d}{Z}\to\mor{d}{Y}$ which is according to
    \cite[Proposition 3]{courbes-rat} an affine bundle and in our situation a
    vector bundle. More precisely define the morphisms ${\rm pr}$ and ${\rm
      ev}$ by the following diagram 
$$\xymatrix{\mor{d}{Y}\times\pu\ar[r]^-{\rm ev}\ar[d]^{\rm pr}&Y\\
\mor{d}{Y}.\\}$$
Then, according to the proof of \cite[Proposition 3]{courbes-rat} the vector
bundle defined by $q$ is associated to the locally free sheaf ${\rm ev}^*{\rm
  pr}_*E$.

Let $f:\pu\to Y$ be an element in $\mor{d}{Y}$ and consider the vector bundle
$Z_f\to \pu$ obtained by pull-back. The associated locally free sheaf is
$f^*E$. The fibre $q^{-1}(f)$ is the set of global sections
$H^0(\pu,f^*E)$. Now if we take $f$ passing through $p(x)$, $p(y)$ and $p(z)$,
then a section $s$ passes through $x,y,z$ if $ev(s)=(e_x,e_y,e_z)$,
where $ev : H^0(\pu,f^*E) \to E_{p(x)} \times E_{p(y)} \times E_{p(z)}$
is the evaluation
at $f^{-1}(p(x)),f^{-1}(p(y)),f^{-1}(p(z))$
and $e_x \in E_{p(x)} , e_y \in E_{p(y)} , e_z \in E_{p(z)}$ are some fixed
elements detemined by $x,y,z$.
By hypothesis there exists a
morphism $f$ such that all factors of $f^*E$ have degree at least 2.
By semi-continuity this is the case for a
general morphism (the space of morphisms $\mor{d}{Y}$ is irreducible, see for
example \cite{courbes-rat}). 
In this case $ev$ is surjective and
$q^{-1}(f)$ is an affine subspace
of $H^0(\pu,f^*E)$. The result follows.
\end{proo}

To prove rationality results with the vector bundle technique, it is therefore 
enough to produce morphisms $f:\pu\to Y$ such that the cohomology group
$H^1(\pu,f^*E\otimes\cO_\pu(-3))$ vanishes. Using the following lemma, 
we will only need to prove this for small degree morphisms.

\begin{lemm}
\label{lemm-degre}
Assume that for a general morphism $f:\p^1 \to Y$ of degree $d$,
the cohomology group $H^1(\pu,f^*E\otimes\cO_\pu(-3))$ vanishes.
Then for any $d'\geq d$, and a general morphism $g:\p^1 \to Y$ of degree $d'$,
the cohomology group $H^1(\pu,g^*E\otimes\cO_\pu(-3))$ vanishes.
\end{lemm}

\begin{proo}
Let $f:\pu\to Y$ a degree $d$ morphism such that 
the cohomology group $H^1(\pu,f^*E\otimes\cO_\pu(-3))$ vanishes.
Denote $y_1$ resp. $y_2,y_3$ the element $f(0)$ resp. $f(1),f(\infty)$.
Let $g$ be a stable morphism $\p^1 \cup \p^1 \to Y$ of degree $d'$
obtained as the map $f$ 
on the first factor and a map $f':\pu\to Y$ of degree $d'-d$ whose image is
not passing through $y_1$, $y_2$ or $y_3$ but meets
$f(\pu)$ in one point 
say $y_4$. For the curve $g(\pu\cup\pu)$ we have the exact sequence 
$0\to f'^*(E(-y_4))\to g^*E\to f^*E\to 0$ and because the points 
$y_1$, $y_2$ and $y_3$ are not on $f'(\pu)$ we get the exact sequence 
$$0\to f'^*(E(-y_4))\to g^*(E(-y_1-y_2-y_3))\to 
f^*(E(-y_1-y_2-y_3))\to 0.$$
Because $E$ is globally generated (see Remark \ref{rema-global-gen}), we get 
the vanishing of $H^1(\pu,f'^*(E(-y_4)))=0$ and by hypothesis we have 
$H^1(\pu,f^*(E(-y_1-y_2-y_3)))=H^1(\pu,f^*E\otimes\cO_\pu(-3))=0$.
We deduce the condition $H^1(\pu\cup\pu,g^*(E(-y_1-y_2-y_3)))=0$. By 
semi-continuity the same holds for a general degree $d'$ stable map 
$g:\pu\to Y$ with three marked points $0$, $1$ and $\infty$ \emph{i.e.} 
we have $H^1(\pu,g^*E\otimes\cO_\pu(-3))=0$. The result follows. Remark that 
we used the irreducibility of the space of space maps in $Y$ (see for 
example \cite{thomsen}).
\end{proo}

\subsection{Proof of the rationality}

In this subsection we prove the rationality of certain Gromov-Witten
varieties using Proposition \ref{prop-fibre} and \ref{prop-rat-fibre}. We
shall deal with the following homogeneous spaces $X$: the Grassmannians
$\G_Q(n,2n)$ of maximal isotropic subspaces in an even dimension vector space
with a non degenerate quadratic form and the Cayley plane $E_6/P_1$. We denote
by $\GW_d(X)$ the Gromov-Witten variety of degree $d$ curves passing through
three general points in $X$. 

\subsubsection{Orthogonal Grassmannians}

Let us start with $X=\G_Q(n,2n)=G/P$ with $P$ the maximal parabolic subgroup 
with $\Sigma(P)=\{\a_n\}$ with notation as in \cite{bou}. Take $Q$ to be the 
maximal parabolic subgroup such that 
$\Sigma(P)\cup i(\Sigma(Q))=\{\a_{n-1},\a_n\}$. The open subset $Z$
of $X$ defined by $Q$ as above has complementary of codimension at least two
since $i(\Sigma(Q))\cap\Sigma(P)$ is empty. Furthermore, by Proposition
\ref{prop-fibre}, the morphism $p:Z\to Y$ defined above is a vector
bundle. The variety $Y$ is easily seen to be isomorphic to $\p^{n-1}$
and the locally free sheaf $E$ associated to $p$ is $\Lambda^2
T_{\p^{n-1}}(-1)$. We conclude by
Proposition \ref{prop-rat-fibre} for $d\geq n-1$ using the following

\begin{prop}
  For $d\geq n-1$, there exists a degree $d$ morphism $f:\pu\to Y$ such that
  the cohomology group $H^1(\pu,f^*E\otimes\cO_\pu(-3))$ vanishes.
\end{prop}

\begin{proo}
By Lemma \ref{lemm-degre}, we only need to construct a degree $n-1$ morphism
$f:\pu \to Y$ with $H^1(\pu,f^*E\otimes\cO_\pu(-3))=0$. But for a general 
degree $n-1$ morphism $f:\pu\to Y$ we have $f^*(T_{\p^{n-1}}(-1))=
\cO_{\pu}(1)^{\oplus n-1}$ and the result follows by taking the second exterior 
power.
\end{proo}

Using the fact that the Gromov-Witten variety of curves passing through 
three points in $\p^n$ is rational and Proposition \ref{prop-rat-fibre}, 
we obtain:

\begin{coro}
  \label{prop-gq-grend-deg}
The Gromov-Witten varieties of degree $d$ curves passing through three general 
points in $G_Q(n,2n)$ is rational for $d\geq n-1$.
\end{coro}

For $d<n/2$, we know by dimension counting, that there is no degree $d$ 
curve passing through three general points in $G_Q(n,2n)$. For $n$ 
even and $d=n/2$, then by a result of A. Buch, A. Kresch and A. Tamvakis 
\cite{BKT}, there is a unique degree $d$ curve passing through 3 general 
points in $G_Q(n,2n)$. We are left to deal with curves of degree $d$ with 
$n/2<d<n-1$. For this we will describe more precisely the geometry of the 
curves through three fixed points. We first make the following remark:

\begin{rema}
\label{kleiman}
Let $F$ be 
a proper closed subset of $\mor{d}{X}$.
Given $x,y,z \in X$ three generic points, denote by $\mor{d}{X,\{x,y,z\}}$
the subscheme of morphisms $f$ such that $f(0)=x,f(1)=y,f(\infty)=z$, and
assume it is not empty.
By Kleiman's transversality Theorem \cite{kleiman}, the scheme
$\mor{d}{X,\{x,y,z\}}$ is equidimensional of dimension $dc_1(X)-2\dim X$
and $\mor{d}{X,\{x,y,z\}} \cap F$ is empty or equidimensional of dimension
$dc_1(X)-2\dim X-\codim\ F$. 
Therefore in any irreducible component of  $\mor{d}{X,\{x,y,z\}}$
there exist elements outside $F$.
\end{rema}

Using this remark and because we want to prove the rationality of the 
Gromov-Witten variety, we will only need to study general degree $d$ 
morphisms. In particular, we may assume that for $f$ such a morphism the 
decomposition of $f^*K$, where $K$ is the tautological subbundle in $X$, is 
the general decomposition. This decomposition is given in the following:

  \begin{lemm}
    Let $K$ be the tautological subbundle on $\G_Q(n,2n)$, then for a general 
morphism $f:\pu\to \G_Q(n,2n)$ of degree $d=n-a$ we have
$f^*K=\cO_\pu(-1)^{\oplus 2a}
\oplus\cO_\pu(-2)^{\oplus n-2a}$.
  \end{lemm}

  \begin{proo}
To prove this result, we decompose the vector space $\C^{2n}$ as a direct sum 
$\oplus_iV_i$ of mutually orthogonal vector spaces such that the restriction 
of the quadratic form is non degenerate on $V_i$ and such that $V_i$ as 
dimension $2n_i$ equal to 4 or 6. We remark that we may define an embedding:
$$\prod_i\G_Q(n_i,2n_i)\to \G_Q(n,2n)$$
by taking the direct sum of the subspaces in each subspace $V_i$. The 
pull-back of the tautological subbundle is the direct sum of the tautological 
subbundles of each factor $\G_Q(n_i,2n_i)$. Now because $\G_Q(n_i,2n_i)$ is 
isomorphic  to $\pu$ or $\p^3$ and the tautological subbundle $K_i$ is 
$\cO_\pu(-1)^{\oplus 2}$ or $\wedge^2(\Omega^1_{\p^3}(1))$, we know the 
decomposition on $\pu$ of the pull-back of a degree $d_i$ morphism 
$f_i:\pu\to \G_Q(n_i,2n_i)$. In particular, we may assume that $f^*K_i$ 
has only factors of the form $\cO_\pu(-1)$ or $\cO_\pu(-2)$ and the result 
follows from the fact that the degree of $K$ on $\G_Q(n,2n)$ is -2.
  \end{proo}

\begin{prop}
\label{prop-gq-petit-deg}
For $n/2<d<n-1$ and $x$, $y$ and $z$ three general points in $\G_Q(n,2n)$, 
the variety of degree $d$ morphisms
$f:\pu\to \G_Q(n,2n)$ with $f(0)=x$, $f(1)=y$ 
and $f(\infty)=z$ is rational.
\end{prop}

\begin{proo}
Let us set $d=n-a$. By the previous Lemma we know that, for a fixed non 
degenerate quadratic form on $\C^{4a}$ (resp. on $\C^{2n-4a}$), there exists 
a map $\cO_\pu(-1)^{2a}\to\cO_\pu\otimes\C^{4a}$ (resp. 
$\cO_\pu(-2)^{n-2a}\to\cO_\pu\otimes\C^{2n-4a}$) defining a degree $a$ 
(resp. $n-2a$) morphism $f_1:\pu\to \G_Q(2a,4a)$ (resp. 
$f_2:\pu\to \G_Q(n-2a,2n-4a)$). Taking $\C^{2n}$ as the orthogonal direct 
sum of $\C^{4a}$ and $\C^{2n-4a}$ we obtain a map 
$\cO_\pu(-1)^{2a}\oplus \cO_\pu(-2)^{n-2a}\to\cO_\pu\otimes\C^{2n}$
defining a degree $n-a$ morphism $f:\pu\to \G_Q(n,2n)$ and such that the 
following diagram commutes
$$\xymatrix{\cO_\pu(-1)^{2a} \ar@{^(->}[r] \ar@{^(->}[d]&
\cO_\pu\otimes\C^{4a}\ar@{^(->}[d]\\
\cO_\pu(-1)^{2a}\oplus \cO_\pu(-2)^{n-2a} \ar@{^(->}[r] \ar@{->>}[d]&
\cO_\pu\otimes\C^{2n} \ar@{->>}[d]\\ 
\cO_\pu(-2)^{n-2a} \ar@{^(->}[r] & \cO_\pu\otimes\C^{2n-4a}.}$$

By Remark \ref{kleiman}, we may assume that a general $f$ in the Gromov-Witten 
variety satisfies the above conditions. In particular, such a map $f$ defines 
a $4a$-dimensional subspace $\C^{4a}$ of $\C^{2n}$ such that the restriction 
of the quadratic form is non degenerate and any subspace $f(t)$ meets 
$\C^{4a}$ 
in dimension $2a$. In particular $\C^{4a}$ meets the three subspaces 
corresponding to $x$, $y$ and $z$ in dimension $2a$ defining therefore three 
points $x'$, $y'$ and $z'$ in $\G_Q(2a,4a)$. The map $f$ induces a 
degree $a$ morphism $f_1:\pu\to \G_Q(2a,4a)$ passing through $x'$, $y'$ and 
$z'$. We already know that there is a unique such morphism. By projection
with respect to $\C^{4a}$ the morphism $f$ defines a degree $n-2a$ morphism 
$f_2:\pu\to \G_Q(n-2a,2n-4a)$ passing through the images, say $x''$, $y''$ 
and $z''$, of $x$, $y$ and $z$ under the projection. By Corollary
\ref{prop-gq-grend-deg} (here the degree is $D=n-2a\geq N-1=n-2a-1$ in 
$G_Q(N,2N)$), we know that the variety of such maps is rational. We are 
therefore left to prove that the variety ${\cal W}$ of $4a$-dimensional 
subspaces $W$ 
of $\C^{2n}$ meeting the three points $x$, $y$ and $z$ in dimension $2a$ is 
rational. 

For this, we need to discuss on the parity of $n$. If $n$ is even, then 
the variety ${\cal W}$ is birational to the Grassmannian $\G(2a,n)$. 
Indeed, if a $4a$-dimensional subspace $W$ is given, its intersection with 
$x$ defines a subspace $x'$ of dimension $2a$ in $x$ (which is of dimension 
$n$). Conversely, if we have a subspace $x'$ of dimension $2a$ in $x$, 
then by projection with respect to $x'$, we send $y$ and $z$ to subspaces 
of dimension $n$ (this is because as $n$ is even, the space $x$, 
 and therefore $x'$, does not meet $y$ and $z$) in $\C^{2n}/x'$ which 
is of dimension $2n-2a$. The images of $y$ and $z$ meet in dimension $2a$ 
and taking the inverse image of this subspace we obtain the desired 
$4a$-dimensional subspace.

If $n$ is odd, we first remark that $P=x\cap y$, $Q=x\cap z$ and $R=y\cap z$ 
are of dimension one.  Denote by $U$ the isotropic
3-dimensional subspace generated 
by $P$, $Q$ and $R$. Let us first prove that  $P$, $Q$ and 
$R$ are contained in $W$ therefore $U\subset W$. Indeed let, as above, be 
$x'$ the $2a$-dimensional intersection $x\cap W$ and consider the projection 
$p_{x'}$ by $x'$. Then $p_{x'}(W)$ is a $2a$-dimensional subspace and 
$p_{x'}(y)$ (resp. $p_{x'}(z)$) is a subspace of dimension $n-1$ or $n$ 
according to the fact that $P$ (resp. $Q$) is contained in $x'$ or not. If we 
define $y'=y\cap W$ and $z'=z\cap W$, then $p_{x'}(y')$ (resp. $p_{x'}(z')$) 
is a subspace of dimension $2a-1$ or $2a$ according to the fact that $P$ 
(resp. $Q$) is contained in $x'$ or not. An easy dimension count shows that 
we have $p_{x'}(y)\cap p_{x'}(z)= p_{x'}(y')\cap p_{x'}(z')$ and is therefore 
contained in $p_{x'}(W)$. In particular we have $p_{x'}(R)\in p_{x'}(W)$, 
therefore $R$ is in $W$. By symmetry the result follows.

%
%The same 
%
%{\bf ** ICI detailler **}
%
%argument as above implies that $P$, $Q$ and 
%$R$ are contained in $W$ therefore $U\subset W$. 

Let us now consider 
$(W\cap U^\perp)/U\subset (\C^{2n}\cap U^\perp)/U$. These spaces are of 
dimension $4a-6$ and $2n-6$ and the spaces $(x\cap U^\perp)/U$, 
$(y\cap U^\perp)/U$ and $(z\cap U^\perp)/U$ are of dimension $n-3$. Therefore 
to determine $(W\cap U^\perp)/U$ we only have to apply the even case.
Now we need to determine $W$ containing 
$W\cap U^\perp$. For this we project with respect to $W\cap U^\perp$. The 
space $W$ becomes a 3-dimensional subspace of $\C^{2n}/W\cap U^\perp$ (which 
is of dimension $2n+3-4a$) and the image of $W$ has to meet the image of $x$, 
of $y$ and of $z$ in dimension 1. These images are the subspaces 
$x/(x\cap W\cap U^\perp)$, $y/(y\cap W\cap U^\perp)$ and 
$z/(z\cap W\cap U^\perp)$ and are of dimension $n+1-2a$. To determine $W$ we 
only need to pick one point in each of the spaces 
$\p(x/(x\cap W\cap U^\perp))$, $\p(y/(y\cap W\cap U^\perp))$ and 
$\p(z/(z\cap W\cap U^\perp))$ and take the inverse image by the projection 
with respect to $W\cap U^\perp$ of the 3-dimensional subspace obtained as the 
span of these 3 points.
\end{proo}

\subsubsection{The Cayley plane}
\label{subsubsection-e6}

In this subsection $X$ is the Cayley plane, that is the variety $G/P$ with 
$G$ of type $E_6$ and $\Sigma(P)=\{\a_1\}$ (again with notation as in 
\cite{bou}). Let us take $Q=P$.
We have $\Sigma(P)\cap i(\Sigma(P))=\emptyset$.
With the notation of Section \ref{section-fibre}, we have a vector bundle 
$p:Z\to Y$. Let us denote by $E$ the locally free sheaf on $Y$ 
corresponding to this vector bundle. 

\begin{lemm}
  We have $Y=\Q^8$ and $E$ is the spinor bundle on $Y$.
\end{lemm}
By a spinor bundle on a smooth quadric $\Q^{2d}$
we mean the vector bundle given by
a spinor representation of $Spin_{2d-2}$, the semi-simple quotient of the
parabolic subgroup of $Spin_{2d}$ stabilising a point in $\Q^{2d}$.
\begin{proo}
The fact that $Y$ is a 8-dimensional quadric follows from the definition of 
$Y$ and the fact that the Levi subgroup $L^{w_0}_P$ of $P^{w_0}$ (resp. 
$L^{w_0}_P\cap P$) is isomorphic to a group of type $D_5$ (resp. to a 
parabolic subgroup of $L^{w_0}_P$ correponding to the root $\a_1$ in $D_5$).

%ICI il faut reprendre pour remplacer $L$ par le groupe quotient reductif.

For the vector bundle $E$, we know that it is given by the representation of 
$L\cap P$ on $U^{w_0}_P/(U^{w_0}_P\cap P)$ where $U^{w_0}_P$ is the unipotent 
radical in $P^{w_0}$. This vector space has for weight under the maximal torus 
of $G$ (of type $E_6$) the highest root of $G$. Therefore it has the weight 
$\varpi_5$ for the group $L^{w_0}_P$ of type $D_5$. By dimension count
$U^{w_0}_P/(U^{w_0}_P\cap P)$ has to be the irreducible representation of 
highest weight $\varpi_5$ and the result follows.
\end{proo}

\begin{lemm}
\label{lemm-spineur}
  Let $Q$ be a smooth quadric and $E \to Q$ a spinor
  bundle.
  For $d\geq 4$, there exists a degree $d$ morphism $f:\pu\to Q$ such that
  the cohomology group $H^1(\pu,f^*E\otimes\cO_\pu(-3))$ vanishes.
\end{lemm}

\begin{proo}
Let us first consider the case when $Q$ is even-dimensional and
denote by $2\delta$ its dimension.
By Lemma \ref{lemm-degre}, we only need to prove this result for $d=4$.
For this we will prove that for a general degree 2 morphism
$f:\pu\to\Q^{2\delta}$, 
we have $f^*E=\cO_\pu(1)^{2^{\delta-1}}$, and take its double. 
Consider in the root system of $Spin_{2\delta}$ 
the two orthogonal roots $\theta$ and $\a_1$ where 
$\theta$ is the highest root. We set $x=x_{\theta}+x_{\a_1}$ where for $\a$ 
root of $Spin_{2\delta}$, the element $x_\a$ is a
non zero element in the root space corresponding to 
$\a$. Take $y=y_{\theta}+y_{\a_1}$ where for $\a$ a root of $Spin_{2\delta}$, 
$y_\a$ is the element with weight $-\a$ and such that 
$[x_\a,y_\a]=\a^\vee$. For $h=[x,y]=\theta^\vee+\a_1^\vee$, 
the triple $(x,h,y)$ is a $\fsl_2$-triple. Denote by
$SL_2(\theta+\a_1)$ the associated  
subgroup of $Spin_{2\delta}$. The orbit of the projectivisation of the highest
weight line under $SL_2(\theta+\a_1)$ is a degree 2 curve and, since
$\scal{\theta^\vee+\a_1^\vee,\a_2} = 0$, all
the weights of the spinor representation evaluated on $\theta^\vee+\a_1^\vee$
are equal to one. The result follows.

For odd-dimensional quadrics, we may either reproduce a similar argument or
use the facts that the restriction of a spinor bundle on $Q^{2\delta}$ to
$Q^{2\delta-1}$ is a spinor bundle on $Q^{2\delta-1}$ and that a rational curve
of degree 4 in a smooth quadric is included in a smooth quadric
of dimension 3.
\end{proo}

By the previous two Lemmas and Proposition \ref{prop-rat-fibre} we get:

\begin{prop}
\label{prop-e6-gd-degre}
For $d \geq 4$ and $x,y$ and $z$ three general points in $E_6/P_1$, the variety
of degree $d$ morphisms $f : \p^1 \to E_6/P_1$ with $f(0)=x,f(1)=y$ and
$f(\infty)=z$ is rational.
\end{prop}

\def \op {{\mathbb {OP}}^1}

We are left with dealing with morphisms of degree $d\leq3$. We prove the 

\begin{prop}
\label{prop-pas-3pts}
  There is no curve of degree 3 passing through 3 general points in $X$.
\end{prop}

\begin{proo}
If such a curve exists, then
the variety of all these curves has dimension 4. In 
particular, by the same proof as Proposition 3.2 in \cite{debarre}, there 
exists a non irreducible degree 3 rational curve through the same three 
general points. This is not possible by the following Lemma
\ref{lemm-reducible-e6}.
\end{proo}

\begin{lemm}
\label{lemm-reducible-e6}
There is no reducible curve of degree 3 in $X$ passing through 3 general 
points.
\end{lemm}

\begin{proo}
Let $x,y,z$ be the three general points.
Up to a permutation, we can assume that there is a curve $C$ of degree 2 
through $x$ and $y$ and a line $l$ through $z$ meeting $C$. To get a
contradiction we
consider the representation $V$ of $E_6$ in which $X$ is embedded
($\dim V = 27$).
Recall that $V$ identifies with the space of order three hermitian
matrices with
octonionic coordinates and $X$ with the variety of matrices with rank one
(see \cite{projective}). We show that for the particular choice
$$
x = \mattt 100 000 000 \ , \ y = \mattt 000 010 000 \ , \ 
z = \mattt 000 000 001 \ \ ,
$$
the existence of $C$ and $l$ leads to a contradiction.

In fact, the union of all conics through $x$ and $y$ is the quadric lying
in the linear subspace given by the matrices of the form $A$ below,
thus $C$ is included in this space. On the other hand $l$ must be
included in the tangent space at $X$ to $z$, which is the linear
subspace of the matrices of the form $B$ below:
$$A=\mattt **0 **0 000,\ B=\mattt 00* 00* ***.$$ 
Thus $l$ and $C$ do not meet.
\end{proo}

We will however prove that there is a quantum to classical $K$-theoretic
principle showing the following results: for a general curve $C$ of degree 3 in
$E_6/P_1$, there is a unique element $\mu \in E_6/P_4$ such that $C$ is included
in the 6-dimensional Schubert variety $S_\mu$ corresponding to $\mu$. There
exists a curve of degree 3 passing through 3 generic points in $S_\mu$
and the space of such curves is rational.

\vskip .2cm

Let us first recall some facts about $X=E_6/P_1$. The following results were 
first proved by \cite{zak} with geometric arguments; they are easy to prove 
once we know that there are three $E_6$-orbits in $X \times X$,
which is proved for example 
in \cite[Proposition 3.16]{CMP}.
First of all a point $\lambda \in E_6/P_6$ defines
a hyperplane section $H_\lambda$ of $X$, singular along a smooth 8-dimensional
quadric that we denote $\op_\lambda$.
There are
three $E_6$-orbits in $X^2$: the pairs $(x,y)$ such that $x=y$, the pairs 
such that 
the line through $x$ and $y$ in $\p V$ is included in $X$
(or equivalently $y \in T_xX$), or
the generic pairs. If $(x,y)$ is a generic pair then there is a unique
$\lambda = \lambda(x,y)$ such that $x,y \in \op_\lambda$.
When $x$ is fixed, the map $y \mapsto \lambda(x,y)$, with image the
smooth 8-dimensional quadric of elements $\lambda \in E_6/P_6$
such that $x \in \op_\lambda$, identifies with the projection
$X \dasharrow \Q^8$ with center $T_xX$.
We denote
$\op(x,y) := \op_{\lambda(x,y)}$.

Similarly there are three $E_6$-orbits in ${(E_6/P_6)}^2$: for generic
$(\lambda,\lambda')$ the intersection $\op_\lambda \cap \op_{\lambda'}$
is transverse and is one point, and for special $(\lambda,\lambda')$,
with $\lambda \not = \lambda'$, this intersection is a maximal isotropic
subspace in both quadrics $\op_\lambda$ and $\op_{\lambda'}$;
thus it is isomorphic to
$\p^4$. 

\vskip .2cm

Let us now explain how a general curve $C$ of degree 3 in $X$
determines an element $\mu(C)$ in $E_6/P_4$. Let $x,y \in C$ be general
points. These are general points in $X$.
Thus they define an element
$\lambda(x,y) \in E_6/P_6$.
In this way we get a rational 
map $\p^2 = S^2\p^1 = S^2C \dasharrow E_6/P_6$.
In fact this rational map has degree one,
as one can see by restricting the lines
in $\p^2=S^2\p^1$ where one element is fixed, say $x=x_0$. In fact,
the elements $\lambda(x_0,y)$ for $y \in C$
are the elements of the image of $C$ under
the projection with center $T_{x_0}X$, which is a curve of degree 1.
It follows that this rational map is in fact a linear embedding
of $\p^2$, and $C$ determines a $\p^2$ in $E_6/P_6$, or an element
$\mu(C) \in E_6/P_4$.

\vskip .2cm

For a given $\mu \in E_6/P_4$,
there are four Schubert varieties homogeneous under the stabiliser
of $\mu$ (a parabolic subgroup conjugated to $P_4$); these
Schubert varieties have dimension
13,10,6,2, and we denote them accordingly 
$S_\mu^{13},S_\mu^{10},S_\mu^6,S_\mu^2$.
Each smaller Schubert variety is the singular locus of the bigger one. 
Remark that $\mu$ determines also a plane in $E_6/P_1$: this is
$S^2_\mu$. 
The variety $S_\mu^{13}$ is the intersection of all the $H_\lambda$
for $\lambda \in \p^2_\mu \subset E_6/P_6$.
Similarly, $S_\mu^{10}$ is the union of the
$\op_\lambda$ for $\lambda \in \p^2_\mu \subset E_6/P_6$.
The variety $S^6_\mu$ is the singular locus of $S_\mu^{10}$.

\begin{lemm}
\label{lemm-inclus}
We have the inclusion $C \subset S_{\mu(C)}^6$.
\end{lemm}
Note that we cannot have $C \subset S_\mu^2$ since otherwise $C$ would only
span a $\p^2$.
\begin{proo}
Let $\mu = \mu(C)$
and let $x \in C$. By definition of $\mu$, if $y \in C$ then
$\lambda(x,y) \in \p^2_\mu$
and therefore
$\op(x,y) \subset S_\mu^{10}$. Thus in particular $x \in S_\mu^{10}$ and
$C \subset S_\mu^{10}$. But let us see that moreover $x$ lies in the singular
locus of $S_\mu^{10}$. In fact since $\op(x,y) \subset S_\mu^{10}$ we deduce
$T_x \op(x,y) \subset T_xS_\mu^{10}$. For $y' \in C$ another element, we
also have $T_x \op(x,y') \subset T_x S_\mu^{10}$. Moreover
$T_x \op(x,y)$ and $T_x \op(x,y')$ meet along $T_x (\op(x,y) \cap \op(x,y'))$
and $\op(x,y) \cap \op(x,y')$ is a linear space of dimension 4,
thus the linear span of $T_x \op(x,y)$ and $T_x \op(x,y')$ has dimension 12.
Thus $x$ is a singular point of $S_\mu^{10}$.
\end{proo}

In the following we denote $S_\mu^6$ simply by $S_\mu$ and we denote
$w = s_2s_6s_5s_4s_3s_1$.
The Picard group of $S_\mu$ has rank 1; however, the group of Weil divisors 
has rank two. Recall that a basis of the group of
Weil divisors is given by the 
Schubert divisors. We denote $D_6 = X(s_6w)$ the singular one and 
$D_2 = X(s_2w)$ the divisor isomorphic to $\p^5$. Therefore,
for a curve in $S_\mu$ 
not meeting the singular locus of $S_\mu$, we may define a bidegree $(a,b)$ by 
intersecting with the Weil divisors $D_6$ and $D_2$. The total degree of the 
curve is $a+b$ since the ample generator of the Picard group is $D_6+D_2$ 
(see for example \cite{small} for more on Cartier and Weil divisors on 
minuscule Schubert varieties).

\begin{lemm}
\label{lemm-courbes-e6}
A generic degree 3 curve $C$ in $X$ does not meet the singular locus 
$S_{\mu(C)}^2$ of $S_{\mu(C)}^6$ and has bidegree $(1,2)$.
The element $\mu(C)$ is the only element $\mu$ such that $C \subset S_\mu$.
\end{lemm}
\begin{proo}
Let us first apply few facts on the space of morphisms from $\pu$ to a 
minuscule Schubert variety as explained in \cite{perrin} to
$S_\mu$. The irreducible
components of the scheme of morphisms from $\pu$ to $S_\mu$
are indexed by the pairs $(a,b)$ of integers such that $a+b=d$: 
indeed a generic curve in the component indexed by $(a,b)$ does not meet the
singular locus of $S_\mu$ and has bidegree $(a,b)$.

Let us now consider the incidence variety $I$ of couples $(C,\mu)$ where 
$\mu \in E_6/P_4$ and $C \subset S_\mu$ is a curve of degree 3.
The projection on the second factor 
$I\to E_6/P_4$ is locally trivial
and has fibers above $\mu$ given by the degree 3 curves in 
$S_\mu$. Thus there are four irreducible components in $I$, according to
the bidegree of $C$ in $S_\mu$; we denote $I_{(a,b)}$ (where $a+b=3$)
these components.
According to \cite{perrin}, the canonical sheaf on
$S_\mu$ is $5D_2+6D_6$ and the space of curves 
of bidegree $(a,b)$ on $S_\mu$ has dimension $5a+6b+6 = 21+b$,
thus $I_{(a,b)}$ has dimension $50+b$.

The set $M$ of couples of the form $(C,\mu(C))$ with $C$ a generic curve of 
degree 3 is the graph of the rational map $\mu$ from $\mor{3}{X}$ to $E_6/P_4$
and thus 
irreducible. It is included in $I$, by Lemma 
\ref{lemm-inclus}. Let
$(a,b)$ be a pair of integers such that $M \subset I_{(a,b)}$.
We first prove that $(a,b) \neq (0,3)$.
Indeed, this would mean that $C$ is contained 
in the divisor $D_2$, which is a $\p^5$. This would imply that the span of $C$ 
is included in $X$, which is not the case for a generic curve $C$.

Since $M \subset I_{(a,b)}$, we have
$52 \leq 50+b$, thus $(a,b) = (1,2)$.
By irreducibility of $I_{(1,2)}$, $M$ is dense in $I_{(1,2)}$,
thus $C \mapsto (C,\mu(C))$ is a rational section of the projection
$p:I_{(1,2)} \to \mor 3X$,
which says that for a generic curve $C$, $\mu(C)$ is the only element $\mu$
such that $C \subset S_\mu$, and implies that $C$ does not meet the singular
locus of $S_{\mu(C)}$.

%We may therefore consider the incidence variety $I'$ of couples $(C,\mu)$ 
%where $\mu \in E_6/P_4$ and $C \subset S_\mu$ is in the component of bidegree 
%$(1,2)$. The projection $I'\to \mor{3}{X}$ is surjective and the map 
%$C\mapsto(C,\mu(C))$ is a rational section of the projection and hence an 
%immersion of an open subset of $\mor{3}{X}$ in $I'$. Since both varieties are 
%irreducible and of the same dimension, we see that the projection $I'\to 
%\mor{3}{X}$ is birational and the last statement follows.

\end{proo}

We now show that the Gromov-Witten variety of $S_\mu$ is rational:

\begin{lemm}
\label{lemm-3pts}
Let $\mu \in E_6/P_4$ and let $x_1,x_2,x_3$ be three generic points in
$S_\mu$. Then the space of curves of bidegree $(1,2)$
in $S_\mu$ which pass through
$x_1,x_2$ and $x_3$ is rational.
\end{lemm}
\begin{proo}
We use the technique developped in Subsection \ref{section-fibre}.
We may assume that $\mu$ is the $B$-stable point in $E_6/P_4$, thus that 
$S_\mu$
admits a $P_4$-action. It follows that $S_\mu$ contains the open subset
$P_4^w/(P_1 \cap P_4^w)$, where $w$ denotes the element
$s_1s_3s_4s_5s_6s_2$. This subset admits a morphism to
$P_4^w/\scal{P_1 \cap P_4^w,R_u(P_4^w)}$, which is also
$L_4^w/(P_1 \cap L_4^w)$.

Now $P_4$ is stable under $s_2,s_6,s_5$, thus we have
$P_4^w = P_4^{s_1s_3s_4}$. The roots of $P_4^w$ resp. $L_4^w$ are thus of 
the form
$$\beta=s_1s_3s_4(\alpha)=s_1s_3s_4\poidsesix abcdef=
\poidsesix{\ d+f-c\ }{\ a+d+f-c\ }
{\ b+d+f-c\ }{\ d\ }{\ e\ }{f},$$ 
with $\alpha$ a root and $c \geq 0$ resp. $c=0$.

Let us compute the roots of $P_1 \cap L_4^w$. In this case we have $c=0$,
and, since $\alpha$ is root, either $f=0$ or $d=0$. Moreover since
$\beta$ must be a root of $P_1$ we have $d+f \geq 0$, thus $d \geq 0$ and
$f \geq 0$. It follows that $L_4^w / (P_1 \cap L_4^w)$ is isomorphic to
$(S_1 \times S_2 \times S_3) / (S_1 \times P_2 \times P_3)$, with
$S_1,S_2$ resp. $S_3$ isomorphic to $SL_3,SL_3$ resp. $SL_2$, and $P_2$
resp. $P_3$ parabolic subgroups of $S_2$ resp. $S_3$ stabilising a line
in $\C^2$ resp. $\C^3$. 
Therefore $L_4^w / (P_1 \cap L_4^w)$ is isomorphic to
$\p^1 \times \p^2$.

Let us now compute the roots in $R_u(P_4^w) / (P_1 \cap R_u(P_4^w))$.
Let $\beta$ be such a root, and write again $\beta = s_1s_3s_4(\alpha)$
with $\alpha$ %=\poidsesix abcdef$ 
a root of $R_u(P_4)$. Since $\alpha$ is a
positive root, we have $d\geq 0$ and $f\geq 0$. Since $\beta$ is not
a root of $P_1$, it is a negative root and thus $d \leq 0$ and $f\leq 0$.
We have $d=f=0$; since $c>0$ ($\alpha \in R_u(P_4)$), it follows that
$c=1$, and $\alpha$ is one of the roots
$$
\poidsesix 001000 \ , \ \poidsesix 011000 \ , \ \poidsesix 111000.
$$

It follows that the fibration
$P_4^w/ (P_1 \cap P_4^w) \to P_4^w / \scal{P_1 \cap P_4^w , R_u(P_4^w)}$
is isomorphic to the total space of the vector bundle
$V_1 \otimes \cO(1,1) \to 
(S_1 \times S_2 \times S_3) / (S_1 \times P_2 \times P_3)$,
where $V_1$ is the 3-dimensional natural representation
of $S_1 \simeq SL_3$. In other words it is the vector bundle
$E:=\cO(1,1)^{\oplus 3}$
on $\p^2 \times \p^1$.

Through 3 points in $\p^2 \times \p^1$ there passes a rational
2-dimensional family of curves
of bidegree $(2,1)$, and the restriction of $E$ to such a curve
is the vector bundle
$\cO(3)^{\oplus 3}$, so that the space of curves passing through 3 points is
rational of dimension 5.
\end{proo}

\subsection{Some more rationality results}

In this subsection, we adress the following
natural question raised in \cite{BM}: is the 
Gromov-Witten variety $\GW_d(X)$ rational for $X$ a homogeneous space and $d$ 
large enough ?
As we have seen, in the cominuscule case, this is always true
since $\GW_d(X)$ is even either empty or rational. 

We quickly explain how, using the previous techniques, one 
can prove such results for some rational homogeneous spaces. We only deal with 
very simple cases where the vector bundles obtained in Proposition 
\ref{prop-fibre} are well understood.

\begin{prop}
\label{prop-ortho}
Let $X$ be an orthogonal Grassmannian, then $\GW_d(X)$ is rational for 
large $d$.
\end{prop}

\begin{proo}
  We start with orthogonal Grassmannians, so let $G$ be semi-simple of type 
$B_n$ or $D_n$ and let $P$ be a maximal parabolic sugroup of $G$. By what we 
already proved, we may assume that $X=G/P$ is not a quadric and therefore 
$X=\G_Q(p,N)$ with $p>1$ and $N=2n+1$ or $2n$. We can choose the parabolic 
$Q$ used in Proposition \ref{prop-fibre} to be maximal and associated to 
the first simple root (with notation as in \cite{bou}). We get an open subset 
$Z$ of $X$ with complement of codimension at least two and a map 
$p:Z\to Y$ where $Y=G'/P'$ with $G'$ of the same type as $G$ with rank one 
less and $P'$ maximal. In other words $Y$ is an orthogonal 
Grassmannian $\G_Q(p-1,N-2)$. The vector bundle $E$ associated 
to $p:Z\to Y$ given by Proposition \ref{prop-fibre} is the tautological 
quotient bundle on $Y$.
 
Then we know that for a general morphism $f:\pu\to Y$ of degree at least 
$N-p-1$, there is no proper vector subspace of the $N-1$ dimensional 
space containing all the subspaces $f(t)$ for $t\in\pu$. This implies that 
$f^*E$ has a decomposition in direct sum of line bundles with positive 
degrees. In particular if we take a double covering $g$ of $f$ we get that 
$g^*E$ has a decomposition in direct sum of line bundles with degrees at least 
2 and we may conclude by induction on the rank of the group and Proposition 
\ref{prop-rat-fibre}. We therefore have the sufficent bound $d\geq 2(N-p-1)$ 
for the rationality of $\GW_d(X)$.
\end{proo}

Recall that a rational homogeneous space with automorphism group $G$ is 
called \emph{adjoint} if it is isomorphic to the closed orbit of the (adjoint) 
action of $G$ on the projectivisation of
its Lie algebra. The adjoint variety for groups of type $A_n$
is isomorphic to the point-hyperplane incidence and has therefore Picard group
of rank two. We will exclude this case. We also exclude the case of groups of 
type $G_2$ for which the method will not work.

\begin{prop}
\label{prop-adj}
  Let $X$ be an adjoint variety for a group $G$ of type different from $A$ and 
$G_2$, then $\GW_d(X)$ is rational for large $d$.
\end{prop}

\begin{proo}
We prove this by case by case analysis. Let us remark that for classical 
groups, the result follows from the previous proposition and the rationality 
results for cominuscule varieties. We therefore only need to deal with the 
exceptional cases.

We first remark, by a simple weight computation, that for all 
adjoint varieties $X$ such that in the Dynkin diagram of ${\rm Aut}(X)$ there 
is a cominuscule vertex (corresponding to a cominuscule root), if we look 
at the open subset $Z$ of $X$ and the vector bundle $p:Z\to Y$ given by 
Proposition \ref{prop-fibre}, then the associated locally free sheaf $E$ on $Y$
is a (non trivial) extension 
$$0\to T_Y^\vee(1)\to E\to \cO_Y(1)\to 0$$ 
where $T_Y$ is the tangent bundle of $Y$. To prove the result we only need to 
prove that there exist morphisms $f:\pu\to Y$ of large degree such that 
$f^*(T_Y^\vee(1))$ has a decomposition as a direct sum of line bundles of 
positive degrees.
For this we only deal with types $E_6$ and $E_7$ as there is no cominuscule 
vertex in the Dynkin diagram of the other exceptional groups. 

In the first case,
the variety $Y$ is a 10-dimensional spinor variety. Its tangent bundle is
$\Lambda^2K^\vee$ where $K$ is the tautological subbundle. If we take a 
general degree 5 morphism $f:\pu\to Y$, then $f^*(K^\vee)=\cO_\pu(2)^5$ thus
$f^*T_Y=\cO_\pu(4)^{10}$ and $f^*(T_Y^\vee(1))=\cO_\pu(1)^{10}$. Therefore we 
may apply Proposition \ref{prop-rat-fibre} as soon as $d\geq 10$ (by taking a 
double covering of $f$ as in the proof of the previous proposition).

In the second case, the variety $Y$ is isomorphic to $E_6/P_1$. If we take a 
general morphism $f:\pu\to Y$ of degree 4, then $f^*T_Y=\cO_\pu(2)^{16}$
(in fact, from the equality $[{\rm pt}]^3=q^4$ in the quantum cohomology, see 
\cite{CMP}, it follows that there is a unique morphism
of degree 4 passing through 4 
general points, which implies the claim).
We get $f^*(T_Y^\vee(1))=\cO_\pu(1)^{16}$. Therefore we may 
apply Proposition 
\ref{prop-rat-fibre} as soon as $d\geq 8$.

For type $E_8$ groups, $E_8/P_8$ has dimension 57 and anticanonical class 29.
We choose for $Q$ the maximal parabolic subgroup 
associated to the simple root $\a_1$. This root is not cominuscule therefore 
the construction explained before Proposition \ref{prop-fibre} gives a tower
of two morphims
$Z \stackrel {p} \to W \stackrel {q} \to Y$,
where $Y$ is isomorphic to a 12-dimensional quadric, $W$ is the total
space of a vector bundle $E$ over $Y$, and $Z$ is
an affine bundle over $W$ whose associated 
direction vector bundle is of the form $q^* F$ with $F$ a vector
bundle on $Y$ (see 
\cite[Proposition 5]{courbes-rat}). These vector bundles are given by the 
ascendent central filtration of the unipotent radical of $Q$. An easy check on 
weights gives us two vector bundles $E$ and $F$ where $F$ is as above 
(therefore of rank 13 and is the tautological quotient bundle on the quadric,
of first Chern class 1) 
and $E$ is a spinor vector bundle (of rank 32 and first Chern class 16).

We can conclude if we know that for a given integer $d$,
and for a generic degree $d$ morphism $f:\p^1 \to Y$,
we have $H^1(\p^1,f^*E) = H^1(\p^1,f^*F) = 0$. For
the vector bundle $F$, we argue as in the proof of Proposition
\ref{prop-ortho}. For $E$, this was in proved in Lemma \ref{lemm-spineur}.

Finally, in the case of $F_4$, we choose the maximal parabolic subgroup
associated to $\a_4$. The argument is very similar to the $E_8$-case: a big
open subset of $F_4/P_1$ (of dimension 15 and anticanonical class 8) is
the total space $Z$ of a tower of affine bundles
$Z \stackrel {p} \to W \stackrel {q} \to Y$. Here $Y$ is a 5-dimensional
quadric, and, with the same notations as above, $F$ is the tautological
quotient bundle on the quadric $Y$ (of rank 6 and first Chern class 1)
and $E$ is the spinor bundle on $Y$ (of rank 4 and first Chern class 2).
Then we conclude as for the type $E_8$-case.
\end{proo}

\section{Quantum to classical principle for equivariant $K$-theory}
\label{setion-k-quant-class}

\subsection{Conditional results of A. Buch and L. Mihalcea}

In this section, we will explain a quantum to classical principle for the 
equivariant quantum $K$-theory of a cominuscule homogeneous space $X$. 
For this, we will follow the ideas of A. Buch and L. Mihalcea \cite{BM} who
proved such results for the Grassmannian varieties and any degree $d$,
and for cominuscule varieties and $d \leq \dmax$,
using a construction in \cite{CMP}. Their technique is valid for any
equivariant quantum $K$-theoretic invariant (\emph{i.e.} $d\geq
d_{\max}$) as soon as the following conjecture is true:

\begin{conj}({\bf Buch-Mihalcea})
  Let $X$ be a cominuscule variety and let $d > d_{\max}$.
  If $x_1,x_2,x_3$ are general points in
  $X$, then the Gromov-Witten variety $GW_d(x_1,x_2,x_3)$ is
  rational.
\end{conj}

We have proved this conjecture in all cases except for one cominuscule
homogeneous space (namely $X=E_6/P_1$) and one degree (namely $d=3$)
for which the variety $GW_d(x_1,x_2,x_3)$ is not rational but empty
(see Proposition  \ref{prop-pas-3pts}). In particular, we will deduce
from the general presentation in \cite{BM} a quantum to  
classical principle for all equivariant quantum $K$-theoretic invariants in 
all these cases except maybe for $E_6/P_1$ and $d=3$.
In Subsection \ref{subsubsection-e6}, we 
also managed to apply their techniques in that case and therefore obtain a
general quantum to classical principle. 

\vskip 0.2 cm

Before stating our result, let us recall the general picture where the 
techniques of A. Buch and L. Mihalcea \cite{BM} may apply. We shall assume 
the following statement about rational curves on the homogeneous space $X$:

\begin{stat}
\label{stat}
Let $d$ be a non negative integer, there exists an element $w_d$ in $W^P$ 
such that
\begin{itemize}
\item for a general degree $d$ morphism $f:\pu\to X$, there exists a
unique translate 
of the Schubert variety $X(w_d)$ containing $f$;
\item for three general points in $X(w_d)$ the variety of degree $d$ morphisms
$f:\pu\to X(w_d)$ passing through these points is rational.
\end{itemize}
\end{stat}

%By the theory on Span and Kernel, A. Buch, A. Kresch and H. Tamvakis proved 
%that this statement is true for any $d$ in a Grassmannian $\G(p,n)$, an 
%orthogonal Grassmannian $\G_Q(n,2n)$ and a Lagrangian Grassmannian 
%$\G_\omega(n,2n)$ and it was proved in \cite{CMP} that this is also the case 
%for $d\leq d_{\max}$ in any cominuscule variety. Furthermore, for all 
%varieties except the Grassmannians and $E_6/P_1$, this statement is true 
%for $d\geq d_{\max}$ with $w_d=w_X$ where $w_X$ is the longest element in 
%$W^P$ as we explain in Section \ref{section-veronese}. This mean that in 
%those varieties, as soon as a degree $d$ morphism passes through two general 
%points, then it also passes through three general points. For Grassmannians, 
%the Span/Kernel techniques prove that the minimal degree for passing through 
%three general points is $\max(p,n-p)$ while $d_{\max}$, the minimal degree 
%for passing through two general points, is $\min(p,n-p)$. However the Span 
%and Kernel techniques prove that the above statement is true. We prove that 
%this statement is also true for $E_6/P_1$, the minimal degree for a morphism 
%to pass to two points is 2 in that case and the minimal degree to pass 
%through three points is 4 which explains the specificity of degree 3
%morphisms.

Once we know that Statement \ref{stat} is true we can, after A. Buch and 
L. Mihalcea, make the following construction. First remark that the 
stabiliser of the Schubert variety $X(w_d)$ is a parabolic subgroup, say 
$Q_d$, and therefore the variety parametrising the translates of the 
Schubert variety $X(w_d)$ is the homogeneous space $Y_d=G/Q_d$. We shall
denote by $Z_d$ the incidence variety between $X$ and $Y_d$. We have the 
following diagram
$$\xymatrix{Z_d\ar[r]^p\ar[d]_q&X\\
Y_d&\\}$$
and for $\mu\in Y_d$, the Schubert variety $S_\mu$ obtained as the translate 
of $X(w_d)$ by $\mu$ is the fiber $q^{-1}(\mu)$. We denote by $Z_d^{(3)}$ 
the triple fiber product $Z_d\times_{Y_d} Z_d\times_{Y_d} Z_d$, we denote by 
$M_d$ the moduli space of degree $d$ stable maps with three marked points 
from a rational curve to $X$ and we define ${\rm B}\ell_d$ by 
$${\rm B}\ell_d=\{(f,S_\mu)\ f\in M_d, \mu\in Y_d\ /\ f\textrm{ is contained 
in }S_\mu\}.$$
As in \cite{BM}, we have the following diagram:
$$\xymatrix{{\rm B}\ell_d\ar[rr]^\pi\ar[d]_\phi&&M_d\ar[d]^{{\rm ev}_i}\\
Z_d^{(3)}\ar[r]^{e_i}&Z_d\ar[r]^p\ar[d]_q&X\\
&Y_d&\\}$$
where the map $\pi$ is just the projection on the first factor, the map $\phi$
forgets $f$ but not the three marked points which become points on $S_\mu$, 
the map ${\rm ev}_i$ for $i=1,2$ or 3 are the three evaluation maps and $e_i$ 
for $i=1,2$ or 3 is the map forgetting the two points in $S_\mu$ of index 
different from $i$.

\begin{theo}
\label{theo-bm}
%  Assume that $\pi$ is birational,
%  that $Z_d^{(3)}$ has rational singularities and that the general
%  fibre of $\phi$ is rational, 
Assume that Statement \ref{stat} holds, then we have the following formula:
$$\chi_{M_d}({\rm ev}_1^*\a\cdot {\rm ev}_2^*\beta\cdot {\rm ev}_3^*\gamma) 
= \chi_{Y_d}(q_*p^*\a\cdot q_*p^*\beta\cdot q_*p^*\gamma) $$
where $\chi_{M_d}:K_T(M_d)\to K_T({\rm pt})$ is the $T$-equivariant Euler 
characteristic and where $\a$, $\beta$ and $\gamma$ are three classes in 
$K_T(X)$.
\end{theo}
\begin{proo}
Remark that if Statement \ref{stat} holds, then $\pi$ is birational, because
of the uniqueness of the translate of $X(w_d)$ containing a general curve.
The variety $Z^{(3)}_d$ is a locally trivial fibration above $Y_d$ with
fibres $X(w_d)^3$; since Schubert varieties have rational singularities,
$Z_d^{(3)}$ therefore also has rational singularities. Finally the second
point of Statement \ref{stat} says that a general fibre of $\phi$ is rational.
Thus we may simply reproduce the proof of Theorem 4.2 in \cite{BM}.
\end{proo}

\subsection{Unconditional results}

Let us now explain how to apply the general framework of the previous
subsection. We have:

\begin{prop}
\label{prop-stat-true}
  Statement \ref{stat} is true for $X$ a cominuscule homogeneous
  space and any degree $d$.
\end{prop}
\begin{proo}
The first point of Statement \ref{stat}
was proved for any cominuscule variety when $d\leq d_{\max}$ in
\cite{BKT,CMP} and one can find in \emph{loc.cit.} a description of the
integers $\dmax$, the Schubert varieties $X(w_d)$, and the varieties $Y_d$.
The remaining results have been proved in this article.
In fact, let us denote by $D_{\max}$ the minimal
degree for a rational curve to pass through three general points of
$X$. We can compute $\dmax$ and $\Dmax$ as
follows:

$$\begin{array}{cccc}
\hline
Type & X &  d_{\max} & D_{\max} \\
\hline
 A_{n-1} & \G(k,n) & \min(p,n-p) & \max(p,n-p) \\
B_n,D_n & \Q^{m} &  2& 2\\
C_n & \G_{\omega}(n,2n) & n+1 & n+1 \\
D_{2n} & \G_Q(2n,4n) & n & n \\
D_{2n+1} & \G_Q(2n+1,4n+2) & n & n+1 \\
E_6 & E_6/P_1 &2 & 4 \\
E_7 & E_7/P_7 & 3 & 3 \\
\hline
\end{array}
$$
For $d\geq D_{\max}$, the varieties $X(w_d)$ and
$Y_d$ are easy to compute: $X(w_d)$
is simply $X$ while $Y_d$ is a point.
In particular the first point of Statement \ref{stat} is trivially true.
The only left cases are
described in the following array (the results come from \cite{BM} for
the Grassmannian variety and Lemmas \ref{lemm-inclus} and
\ref{lemm-3pts} for the $E_6/P_1$) where we assume that
$d\in[d_{\max}+1,D_{\max}-1]$:
$$\begin{array}{cccccccc}
\hline
Type & X & d & X(w_d) & Y_d \\
\hline
 A_{n-1} & \G(k,n) & \dmax < d < \Dmax & \G(\dmax,d) & \G(d+d_{\max},n) \\
E_6 & E_6/P_1 & d=3 & X(s_6s_5s_4s_3s_1s_2) & E_6/P_4 \\
\hline
\end{array}
$$
The first point of Statement \ref{stat} is proved in \cite{BKT} in the case
of Grassmannians and is contained in Lemma \ref{lemm-courbes-e6} in the
case of $E_6/P_1$.

For the second point of Statement \ref{stat}, namely the rationality
of the Gromov-Witten varieties, we observe that when $d \leq \dmax$ this
variety is one point, by \cite{BKT,CMP}, so assume $d > \dmax$.
The case of Grassmannians of type A
was done in \cite{BM}. The case of spinor varieties is done by Corollary
\ref{prop-gq-grend-deg} and Proposition \ref{prop-gq-petit-deg}. The case
of symplectic Grassmannians is covered by Theorem \ref{theo-veronese}. For
$E_6/P_1$, we use Lemma \ref{lemm-3pts} for $d=3$ and Lemma \ref{lemm-spineur}
for $d \geq 4$.
\end{proo}

\noindent
By Theorem \ref{theo-bm}, we thus get:

\begin{theo}
For $X$ a cominuscule homogeneous space and $d$ any non negative
integer, we have the following formula:
$$\chi_{M_d}({\rm ev}_1^*\a\cdot {\rm ev}_2^*\beta\cdot {\rm ev}_3^*\gamma) 
= \chi_{Y_d}(q_*p^*\a\cdot q_*p^*\beta\cdot q_*p^*\gamma) $$
where $\chi_{M_d}:K_T(M_d)\to K_T({\rm pt})$ is the $T$-equivariant Euler 
characteristic and where $\a$, $\beta$ and $\gamma$ are three classes in 
$K_T(X)$.
\end{theo}

As by definition, the equivariant quantum $K$-theoreric Gromov-Witten 
invariant $I_d^T(\a,\beta,\gamma)$ is $\chi_{M_d}({\rm ev}_1^*\a\cdot 
{\rm ev}_2^*\beta\cdot {\rm ev}_3^*\gamma)$, we deduce the equivariant
$K$-quantum to classical principle:

\begin{coro}
\label{coro-principe}
Let $X$ be a cominuscule homogeneous space and $d$ any integer.
We have the formula \hskip 0.1 cm $I_d^T(\a,\beta,\gamma)= 
\chi_{Y_d}(q_*p^*\a\cdot q_*p^*\beta\cdot q_*p^*\gamma).$
\end{coro}

For $d\geq D_{\max}$, since $Y_d$ is a point, $Z_d^{(3)}$ is equal to
$X \times X \times X$, and we get (compare with Consequence 
7.2 in \cite{BM}): 

\begin{coro}
For $d\geq D_{\max}$, we have the formula \hskip 0.1 cm $I_d^T(\a,\beta,
\gamma)= \chi_{X}(\a)\cdot\chi_{X}(\beta)\cdot \chi_{X}(\gamma).$
\end{coro}

\bigskip\noindent
Pierre-Emmanuel {\sc Chaput}, \\
{\it Laboratoire de Math{\'e}matiques Jean Leray,} 
UMR 6629 du CNRS, UFR Sciences et Techniques,  2 rue de la
Houssini{\`e}re, BP 92208, 44322 Nantes cedex 03, France and \\
{\it Hausdorff Center for Mathematics,}
Universit{\"a}t Bonn, Villa Maria, Endenicher
Allee 62, 
53115 Bonn, Germany.

\noindent {\it email}: \texttt{pierre-emmanuel.chaput@math.univ-nantes.fr}.

\medskip\noindent
Nicolas {\sc Perrin}, \\
{\it Hausdorff Center for Mathematics,}
Universit{\"a}t Bonn, Villa Maria, Endenicher
Allee 62, 
53115 Bonn, Germany and \\
{\it Institut de Math{\'e}matiques de Jussieu,} 
Universit{\'e} Pierre et Marie Curie, Case 247, 4 place
Jussieu, 75252 Paris Cedex 05, France.

\noindent {\it email}: \texttt{nicolas.perrin@hcm.uni-bonn.de}.

\end{document}

%%% Local Variables: 
%%% mode: latex
%%% TeX-master: t
%%% End: